\documentclass[journal]{IEEEtran}

\usepackage{cite}
\usepackage{amsmath,amssymb,amsfonts}
\usepackage{algorithm,algorithmic}
\usepackage{cases}
\usepackage{graphicx}
\usepackage{textcomp}
\usepackage{dsfont}
\usepackage{enumerate}
\usepackage{url}

\usepackage{theorem}
\theorembodyfont{\rmfamily}

\newtheorem{prob}{Problem}
\newtheorem{lem}{Lemma}
\newtheorem{defin}{Definition}
\newtheorem{theorem}{Theorem}
\newtheorem{prop}{Proposition}

\newcommand{\mcX}{\mathcal{X}}
\newcommand{\mcA}{\mathcal{A}}
\newcommand{\mc}[1]{\mathcal{#1}}
\newcommand{\Xa}{\mathcal{X}_{\rm a}}
\newcommand{\bR}{\mathbb{R}}

\newcommand{\Pih}{\Pi^{\rm h}}
\newcommand{\Prob}{\mathbb{P}}
\newcommand{\Exp}{\mathbb{E}}
\newcommand{\bff}{f}
\newcommand{\hPih}{\hat{\Pi}^{\rm h}}
\newcommand{\hpi}{\hat{\pi}}
\newcommand{\hPim}{\hat{\Pi}^{\rm m}}
\newcommand{\ul}[1]{\underline{#1}}

\usepackage[caption=false,font=footnotesize]{subfig}

\newcommand{\hs}{& \hspace{-2mm}}

\begin{document}
\title{Attack Impact Evaluation for Stochastic Control Systems through Alarm Flag State Augmentation}
\author{Hampei Sasahara, \IEEEmembership{Member, IEEE}, Takashi Tanaka, \IEEEmembership{Member, IEEE}, and Henrik Sandberg, \IEEEmembership{Senior Member, IEEE}
\thanks{H. Sasahara is with the Department of Systems and Control Engineering, School of Engineering, Tokyo Institute of Technology, Tokyo 152-8552, Japan (e-mail: sasahara@sc.e.titech.ac.jp). }
\thanks{T. Tanaka is with the Department of Aerospace Engineering and Engineering Mechanics, Cockrell School of Engineering, The University of Texas at Austin, TX 78712, USA (e-mail: ttanaka@utexas.edu).}
\thanks{H. Sandberg is with the Division of Decision and Control Systems, School of Electrical Engineering and Computer Science, KTH Royal Institute of Technology, Stockholm SE-100 44, Sweden (e-mail: hsan@kth.se).}
\thanks{This work was supported in part by JSPS KAKENHI Grant Number 22K21272, Swedish Research Council grant 2016-00861, and Swedish Contingencies Agency.}}

\maketitle

\begin{abstract}
This note addresses the problem of evaluating the impact of an attack on discrete-time nonlinear stochastic control systems.
The problem is formulated as an optimal control problem with a joint chance constraint that forces the adversary to avoid detection throughout a given time period.
Due to the joint constraint, the optimal control policy depends not only on the current state, but also on the entire history, leading to an explosion of the search space and making the problem generally intractable.
However, we discover that the current state and whether an alarm has been triggered, or not, is sufficient for specifying the optimal decision at each time step.
This information, which we refer to as the alarm flag, can be added to the state space to create an equivalent optimal control problem that can be solved with existing numerical approaches using a Markov policy.
Additionally, we note that the formulation results in a policy that does not avoid detection once an alarm has been triggered.
We extend the formulation to handle multi-alarm avoidance policies for more reasonable attack impact evaluations, and show that the idea of augmenting the state space with an alarm flag is valid in this extended formulation as well.
\end{abstract}

\begin{IEEEkeywords}
Attack impact evaluation, chance constraint, control system security, stochastic optimal control.
\end{IEEEkeywords}

\section{INTRODUCTION}
\label{sec:introduction}

The security of control systems has become a pressing concern due to the increased connectivity of modern systems.
There have been indeed numerous reported incidents in industrial control systems~\cite{History2018Hemsley}, and some critical infrastructures that have been seriously damaged~\cite{Nicolas2011Stuxnet,CISA2014,CISA2017,CISA2018}.
Security risk assessment is a crucial step in preventing such incidents from occurring.
For general information systems, risk assessment is typically conducted by identifying potential scenarios, quantifying their likelihoods, and evaluating the potential impacts~\cite{Kaplan1981On,Sridhar2012Cyber}.

Evaluating the impact of attacks to control systems is a challenging task because the defender must specify harmful attack signals.
Many studies have attempted to address this issue by quantifying impact as the solution of a constrained optimal control problem or the reachable set under a constraint (see e.g.,~\cite{Teixeira2015Secure,Murguia2018Reachable}).
These constraints often describe the stealthiness of the attack throughout the considered time horizon.
However, a common problem with these formulations is that they are limited in the types of systems that they can handle.
Existing studies often assume specific forms of the attack detector, such as the $\chi^2$ detector and the cumulative sum (CUSUM) detector~\cite{Murguia2019On}.
However, other types of detectors are also used in practice with distinct properties in terms of detection performance and computational efficiency~\cite{Alvaro2011Attacks}.
While some works provide a universal bound for all possible detectors by using the Kullback-Leibler divergence between the observed output and the nominal signal~\cite{Bai2017Data}, this approach can lead to overly conservative evaluations when the implemented detector is specified.


This study aims to provide a framework for evaluating the impact of attacks that can be applied to a wide range of systems and detectors.
We use a constrained optimal control formulation, in which the stealth condition is represented as a temporally joint chance constraint.
This constraint limits the probability that an alarm is triggered at least once throughout the entire time horizon.
However, because the chance constraint is joint over time, the optimal policy depends not only on the current state, but also the entire history.
As a result, the size of the search space increases exponentially with the length of the horizon, making the problem intractable even for small instances.

In this note, we propose a reformulation of the attack impact evaluation problem in a computationally tractable form.
Our key insight is that the information of whether an alarm has been triggered, in addition to the current state, is sufficient to identify the worst attack at each time step.
We refer to this binary extra information as an \emph{alarm flag}.
By augmenting the alarm flag state to the original state space, we show that the optimal value can be attained by Markov policies.
The reformulated problem is a standard constrained stochastic optimal control problem, which can be solved using exiting numerical methods if the dimension of the spaces is not too large.
Additionally, we note that the adversary does not avoid detection once an alarm has been triggered at least once.
However, this behavior may not be reasonable in practice due to the presence of false alarms.
To address this, we generalize the formulation to handle multi-alarm avoidance policies, providing a more realistic attack impact evaluation.
We also demonstrate that the idea of flag state augmentation is valid in this extended formulation.

\subsection*{Related Work}

The attack impact evaluation problem for control systems has been considerably studied~\cite{Teixeira2015Secure,Mo2016On,Bai2017Data,Umsonst2017Security,Hirzallah2018Computation,Murguia2018Reachable,Chen2018Optimal,Murguia2019On,Teixeira2019Optimal,Milosevic2019Estimating,Fang2019Stealthy,
Sui2020The,Wang2022Optimal,Khazraei2022Resiliency}.
These works formulate the problem as a constrained optimal control problem, but the computation approaches differ based on the type of system, detector, and the objective function used to quantify the attack impact.
To the best of our knowledge, this study is the first work to handle general systems and detectors.

Our idea of alarm flag state augmentation is to add information sufficient for determining the optimal decision at each time step.
A similar concept has been proposed in previous studies especially in the context of risk-averse Markov decision process (MDP)~\cite{Bauerle2011Markov,Haskell2015Convex,Chow2017Risk}.
The work~\cite{Bauerle2011Markov} treats a non-standard MDP where the objective function is given not by expectation but by conditional-value-at-risk (CVaR), to which value iteration can be applied by considering an augmented state space for CVaR.
In~\cite{Haskell2015Convex}, this idea is generalized to chance-constrained MDP.
The work~\cite{Chow2017Risk} proposes risk-aware reinforcement learning based on state space augmentation.
Moreover, linear temporal logic specification techniques can handle general properties, such as safety, invariant, and liveness for discrete event systems~\cite{Baier2008Principles}.
However, our study provides a clear interpretation of the augmented state in the context of control system security, leading to a reasonable extension to the multi-alarm avoidance problem discussed in Sec.~\ref{sec:multi}.
Additionally, we consider a continuous state space, whereas existing studies mainly focus on finite or discrete state spaces.

Temporally joint chance constraints in optimal control have also been studied~\cite{Ono2012Joint,Ono2015Chance,Thorpe2022Data},
but these methods rely on approximating the chance constraint.
Furthermore, they do not discuss the state space augmentation of the decision at each time step.
Finally, a continuous-time optimal control problem with a joint-chance constraint is considered in~\cite{Patil2022ChanceA,Patil2022ChanceB} although the process stops once the state reaches the unsafe region in their formulation.

A preliminary version of this note has been presented in~\cite{Sasahara2022Attack}, but it only considers finite MDPs, and the current study extends it to continuous spaces.
Furthermore, the previous study did not included the theoretical claims and the proofs of Propositions~\ref{prop:1} and~\ref{prop:2} in this note.

\subsection*{Organization and Notation}
This note is organized as follows.
Sec.~\ref{sec:prob} defines the system model, clarifies the threat model, and formulates the attack impact evaluation problem.
In Sec.~\ref{sec:approach}, the difficulty of the formulated problem is explained.
Subsequently, we propose a problem reformulation in a tractable form by introducing the alarm flag state augmentation.
Sec.~\ref{sec:multi} first provides a characterization of the optimal policy after an alarm is triggered.
Based on the observation, we extend the formulation to be able to handle multi-alarm avoidance policies and show that the proposed idea is still valid in the extended variant.
In Sec.~\ref{sec:num}, the theoretical results are verified through numerical simulation, and finally, Sec.~\ref{sec:conc} concludes and summarizes this note.

We denote the set of real numbers by $\mathbb{R}$,
the $n$-dimensional Euclidean space by $\mathbb{R}^n$,
the $t$-ary Cartesian power of the set $\mcX$ for a positive integer $t$ by $\mcX^t$,
the complement of a set $\mc{X}$ by $\mc{X}^{\rm c}$,
the tuple $(x_0,\ldots,x_t)$ by $x_{0:t}$,
and the Borel algebra of a topological space $\mc{X}$ by $\mc{B}_{\mc{X}}$.

\section{ATTACK IMPACT EVALUATION PROBLEM}
\label{sec:prob}
\subsection{System Model}
\label{subsec:sys}

Consider a discrete-time nonlinear stochastic control system of the form
\[
 x_{t+1}=f(x_t,a_t,w_t),\quad t=0,1,\ldots
\]
with $f:\mcX\times\mcA\times\mc{W}\to\mcX$ where $x_t\in\mcX \subset \bR^n$ is the state, $a_t\in\mcA\subset \bR^m$ is the attack signal, and $w_t\in\mc{W}\subset \bR^\ell$ is an independent random process noise.
The distribution of the initial state is denoted by $p_0(dx_0)$.
An attack detector equipped with the control system triggers an alarm when the state reaches the alarm region $\Xa\subset\mcX$.

\begin{figure}[t]
  \centering
  \includegraphics[width=0.98\linewidth]{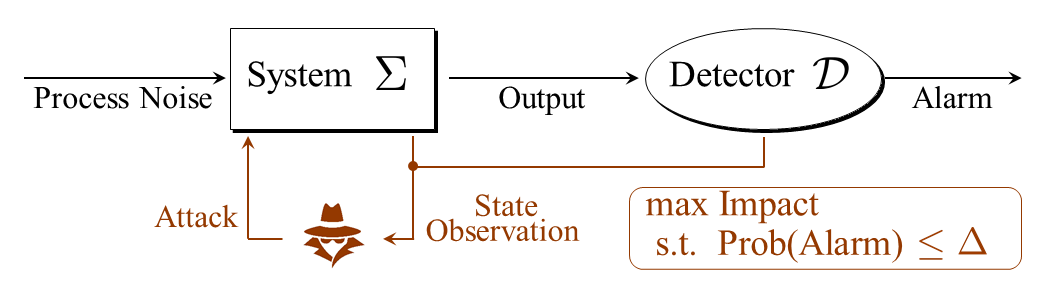}
  \caption{
  Block diagram of a stochastic system with a detector that has been compromised by an attacker.
  The attacker aims at maximizing the damage caused by the attack signal while avoiding triggering the alarm.
  }
  \label{fig:example}
\end{figure}

{\it Remark:} This model includes control systems with the typical cascade structure illustrated by Fig.~\ref{fig:example}.
Let the dynamics of the control system $\Sigma$ and the detector $\mc{D}$ be given by
\[
 \Sigma: z_{t+1}=f^{\rm s}(z_t,a_t,w_t),\quad \mc{D}:\left\{
 \begin{array}{ll}
  z^{\rm d}_{t+1}\hs =f^{\rm d}(z_t,z^{\rm d}_t),\\
  \delta_t\hs =h(z_t,z^{\rm d}_t),
 \end{array}
 \right.
\]
respectively, with the state spaces $\mc{Z}$ and $\mc{Z}^{\rm d}$.
The binary signal $\delta_t\in\{0,1\}$ describes whether an alarm is triggered, or not, at the time step.
It is clear that the cascade system can be described in the form by taking $\mc{X}:=\mc{Z}\times\mc{Z}^{\rm d}$ and $\Xa:= h^{-1}(\{1\})$.

\subsection{Threat Model}
In this study, we consider the following threat model:

\begin{itemize}
\item The adversary has succeeded in intruding into the system and can execute any attack signal in a probabilistic manner at every time step.
\item The adversary has perfect model knowledge.
\item The adversary possesses infinite memory and computation resources.
\item The adversary can observe the state at every time step.
\item The attack begins at $t=0$ and ends at $t=T-1$.
\end{itemize}

The threat model implies that the adversary can implement an arbitrary history-dependent randomized policy $\pi\in\Pih$, where $\pi=(\pi_t)_{t=0}^{T-1}$ is a tuple of policies at every time step.
Let $h_t\in\mc{H}_t:=(\mcX\times\mcA)^{t-1}\times \mc{X}$ denote the history.
The policy at each time step is a stochastic kernel on $\mcA$ given $\mc{H}_t$ denoted by $\pi_t(da_t|x_{0:t},a_{0:t-1})$.
It is well known that a policy $\pi$ uniquely induces a probability measure $\Prob^\pi$ on $\mc{H}_T$~\cite[Chap.~11]{Hinderer1970Foundations}.
The probabilities are specifically given by
\[
\begin{array}{l}
 \Prob^\pi(X_0,A_0,X_1,A_1,\ldots,X_T) := \int_{X_0} \int_{A_0} \cdots \int_{X_T}\\
 \quad \quad  \quad p_T(dx_T|x_{T-1},a_{T-1})
  \cdots \pi_0(da_0|x_0) p_0(dx_0),
\end{array}
\]
where $X_t\in \mc{B}_{\mcX}$ and $A_t\in \mc{B}_{\mcA}$ and $p_t(dx_t|x_{t-1},a_{t-1})$ is the state transition kernel on $\mcX$ given $\mcX\times\mcA$ induced by $f$ and the distribution of $w_{t-1}$~\cite[Chap.~8]{Bertsekas1996Stochastic}.
The expectation operator with respect to $\Prob^\pi$ is denoted by $\Exp^\pi$.

The objective of the adversary is to maximize a cumulative attack impact while avoiding detection.
Let the impact be quantified as
\[
\textstyle{
 J(\pi):=\Exp^\pi\left[ \sum_{t=0}^{T-1} g_t(x_t,a_t) + g_T(x_T) \right].
 }
\]
The adversary keeps herself stealthy over the considered time period $\mc{T}:=\{0,\ldots,T\}$.
Specifically, the probability the event that an alarm is triggered at some time step defined by
\[
 (\lor_{t\in\mc{T}} X_t \in \Xa):=\{(x_{0:T},a_{0:T-1}): \exists t\in\mc{T}\ {\rm s.t.}\ x_t\in\Xa \}
\]
is made less than or equal to a given constant $\Delta \geq 0$.

\subsection{Problem Formulation}

The attack impact evaluation problem is formulated as a stochastic optimal control problem with a temporally joint chance constraint:
\begin{prob}\label{prob:single}
The attack impact evaluation problem is given by
\begin{equation}\label{eq:prob_ori}
\begin{array}{cl}
\displaystyle{\max_{\pi\in\Pih}} & J(\pi)\\
 {\rm s.t.} & \Prob^\pi( \lor_{t\in\mc{T}} X_t \in \Xa) \leq \Delta
\end{array}
\end{equation}
with a given constant $\Delta\geq 0$.
\end{prob}
In the subsequent section, we explain its difficulty and propose an equivalent reformulation in a tractable form.

\section{EQUIVALENT REFORMULATION to TRACTABLE PROBLEM}
\label{sec:approach}

\subsection{Alarm Flag State Augmentation}

It is well known that Markov policies, which depend only on the current state, can attain the optimal value for unconstrained stochastic optimal control problems~\cite[Proposition 8.1]{Bertsekas1996Stochastic}.
However, Problem~\ref{prob:single} has a temporally joint chance constraint that cannot be decomposed with respect to time steps.
Hence Markov policies cannot attain the optimal value in general, an example of which is provided in the Appendix. 
As a result, the size of the search space grows exponentially with the time horizon length,
making the problem intractable even for small instances.

The key idea in this paper to overcome this challenge is to augment the alarm history information into the state space.
We define the augmented state space and the induced augmented system next.
\begin{defin}
The augmented state space of $\mcX$ for Problem~\ref{prob:single} is defined as
\[
 \hat{\mcX}:=\mcX\times \mc{F},\quad
 \mc{F}:=\{0,1\}.
\]
The augmented system is defined as
\[
\left\{
\begin{array}{cl}
 x_{t+1} \hs = f(x_t,a_t,w_t),\\
 \bff_{t+1} \hs = \left\{
 \begin{array}{cl}
 1 & {\rm if}\ x_{t+1}\in\Xa\ {\rm or}\ \bff_t=1,\\
 0 & {\rm otherwise}.
 \end{array}
 \right.
\end{array}
\right.
\]
\end{defin} 

The augmented state $\bff_t\in\mc{F}$ is referred to as the alarm flag, since $\bff_t=1$ indicates that the alarm has been triggered before the time step $t\in\mc{T}$, whereas $\bff_t=0$ indicates otherwise.
For the augmented system, we denote the set of histories by $\hat{\mc{H}}_t:=(\mcX\times\mc{F}\times\mcA)^{t-1}\times (\mcX\times\mc{F})$,
the set of history-dependent randomized policies by $\hPih$,
the probability measure on $\hat{\mc{H}}_T$ induced by $\hpi\in\hPih$ by $\Prob^{\hpi}$,
and the expectation operator with respect to $\Prob^{\hpi}$ by $\Exp^{\hpi}$.

By using the alarm flag, we can rewrite the temporally joint chance constraint in~\eqref{eq:prob_ori} as an isolated chance constraint on the state only at the final time step.
It is intuitively true that $(\lor_{t\in\mc{T}} X_t \in \Xa)$, the event that an alarm is triggered at some time step, is equivalent to $\bff_T^{-1}(\{1\})$, the event that the alarm flag takes the value $1$ at the final time step.
This idea yields the reformulated problem
\begin{equation}\label{eq:prob_r1}
\begin{array}{cl}
\displaystyle{\max_{\hpi\in\hPih}} & J(\hpi)\\
 {\rm s.t.} & \Prob^{\hpi}(\bff_T=1) \leq \Delta.
\end{array}
\end{equation}

The most significant aspect of this formulation is that the chance constraint depends on the marginal distribution with respect to the final time step only.
Hence, the optimal value of~\eqref{eq:prob_r1} can be attained by Markov policies for augmented state space $\hat{\mcX}$~\cite[Proposition 8.1]{Bertsekas1996Stochastic}.
Thus, the problem~\eqref{eq:prob_r1} can be reduced to
\begin{equation}\label{eq:prob_r2}
\begin{array}{cl}
\displaystyle{\max_{\hpi\in\hPim}} & J(\hpi)\\
 {\rm s.t.} & \Prob^{\hpi}(\bff_T=1) \leq \Delta.
\end{array}
\end{equation}
where the search space is replaced with $\hPim$, the set of Markov policies for the augmented system.

\subsection{Equivalence}

We justify the reformulation in a formal manner.
First, we show the following lemma.
\begin{lem}\label{lem:1}
For any $\hpi\in\hPih,$ there exists $\pi\in\Pih$ such that
\[
 \Prob^{\hpi}(X_{0:t},A_{0:t})=\Prob^\pi(X_{0:t},A_{0:t}),\quad \forall t\in \mc{T}
\]
for any $X_{0:t}\in\mc{B}_{\mcX^{t+1}}$ and $A_{0:t}\in\mc{B}_{\mcA^{t+1}}$.
\end{lem}
\begin{IEEEproof}
We say $\bff_{0:t}\in\mc{F}^{t+1}$ to be consistent with $x_{0:t}\in\mcX^{t+1}$ when $\bff_{0:t}$ satisfies
\[
 \left\{
 \begin{array}{ll}
 \bff_t = 1, & \forall t\geq \ul{t}:=\min\{\tau\in\{0,\ldots,t\}:\bff_\tau=1\},\\
 x_t \notin \Xa, & \forall t < \ul{t},\\
 x_{\ul{t}} \in \Xa.
 \end{array}
 \right.
\]
It is clear that a state trajectory $x_{0:t}$ deterministically specifies the consistent alarm flag trajectory $\bff_{0:t}$, denoted by $\bff^{\rm c}_{0:t}(x_{0:t})$.
Note that $P^{\hpi}_{\bff}(\bff_t|x_{0:t})=1$ if $\bff_t=\bff^{\rm c}_{t}(x_{0:t})$ and zero otherwise, where $P^{\hpi}_{\bff}$ denotes the probability mass function on $\bff_{t}$ conditioned on $x_{0:t}$ under the policy $\hpi$.

For a given $\hpi\in\hPih$, determine $\pi\in\Pih$ by
\begin{equation}\label{eq:pi_cons}
 \pi_t(da_t|x_{0:t},a_{0:t-1}):=\hpi_t(da_t|x_{0:t},\bff^{\rm c}_{0:t}(x_{0:t}),a_{0:t-1})
\end{equation}
for $t=0,\ldots,T-1$.
We confirm next that the policy above satisfies the condition in the lemma statement.
From the definition of $\Prob^{\hpi}$ and~\eqref{eq:pi_cons}, we have
\[
\begin{array}{l}
 \Prob^{\hpi}(X_{0:t},A_{0:t})\\
 = \int_{X_0}\int_{A_0}\sum_{\bff_0\in\mc{F}} \cdots \int_{X_t}\int_{A_t}\sum_{\bff_t\in\mc{F}}\\
 \quad \quad P_{\bff}(\bff_t|x_{0:t})\hpi_t(da_t|x_{0:t},\bff_{0:t},a_{0:t-1})p_t(dx_t|x_{t-1},a_{t-1})\\
 \quad \quad \cdots P_{\bff}(\bff_0|x_0)\hpi_0(da_0|x_0,\bff_0) p_0(dx_0)\\
 =\int_{X_0}\int_{A_0}\cdots \int_{X_t}\int_{A_t}\\
 \quad \quad \hpi_t(da_t|x_{0:t},\bff^{\rm c}_{0:t}(x_{0:t}),a_{0:t-1})p_t(dx_t|x_{t-1},a_{t-1})\\
 \quad \quad \cdots \hpi_0(da_0|x_0,\bff^{\rm c}_0(x_0)) p_0(dx_0)\\
 = \int_{X_0}\int_{A_0}\cdots \int_{X_t}\int_{A_t}\\
 \quad \quad \pi_t(da_t|x_{0:t},a_{0:t-1})p_t(dx_t|x_{t-1},a_{t-1})\\
 \quad \quad \cdots \pi_0(da_0|x_0) p_0(dx_0)\\
 = \Prob^{\pi}(X_{0:t},A_{0:t}).
 \end{array}
\]
\end{IEEEproof}

Lemma~\ref{lem:1} implies that the stochastic behaviors of the original system and the augmented one are identical with appropriate policies related to each other through~\eqref{eq:pi_cons}.

The following theorem is the main result of this paper.

\begin{theorem}\label{thm:1}
The optimal values of the problems~\eqref{eq:prob_ori},~\eqref{eq:prob_r1}, and~\eqref{eq:prob_r2} are equal.
\end{theorem}
\begin{IEEEproof}
Denote the optimal values of~\eqref{eq:prob_ori},~\eqref{eq:prob_r1}, and~\eqref{eq:prob_r2} by $J^\ast, \hat{J}^\ast,$ and $\hat{J}^{{\rm m}\ast}$, respectively.
We first show $J^{\ast}=\hat{J}^{\ast}$.
Since the policy set of of the augmented system includes that of the original system, $J^\ast\leq \hat{J}^\ast$ clearly holds.
Fix a feasible policy $\hpi\in\hPih$ for~\eqref{eq:prob_r1} and take the corresponding policy $\pi\in\Pih$ for the original system according to~\eqref{eq:pi_cons}.
From Lemma~\ref{lem:1}, the marginal distributions of the state and the action with the policies coincide.
From the dynamics of $\bff_t,$ we have $\Prob^{\hpi}(\bff_T=1)=\Prob^\pi(\lor_{t\in\mc{T}} X_t \in \Xa)\leq \Delta$, and hence $\pi$ is feasible in~\eqref{eq:prob_ori}.
Therefore, $J^\ast\geq \hat{J}^{\ast},$ which leads to $J^\ast=\hat{J}^\ast$.
Finally, $\hat{J}^\ast=\hat{J}^{{\rm m}\ast}$ is a direct conclusion of~\cite[Proposition 8.1]{Bertsekas1996Stochastic}.
\end{IEEEproof}

Theorem~\ref{thm:1} justifies the reformulation from~\eqref{eq:prob_ori} to~\eqref{eq:prob_r2}.
The reformulated problem~\eqref{eq:prob_r2} can be solved using existing numerical methods if the dimension of the spaces is not too large~\cite{Munos2002Variable}.

\section{EXTENSION: MULTI-ALARM AVOIDANCE POLICY}
\label{sec:multi}

In this section, we observe that the adversary does not avoid detection once an alarm has been triggered at least once based on the previous section's result.
However, this behavior may not be reasonable because of the presence of false alarms.
We generalize the formulation to be able to handle multi-alarm avoidance policies, providing a more reasonable evaluation of the attack impact.

\subsection{Optimal Policy after Alarm Triggered}
\label{subsec:opt_pol_after_alarm}

We observe that the optimal policy after an alarm is triggered is characterized using an optimal policy of an unconstrained problem.
Consider the problem
\begin{equation}\label{eq:prob_unconstrained}
\max_{\pi\in\Pi^{\rm m}} J(\pi)
\end{equation}
and assume that there exists a unique optimal Markov policy, denoted by $\pi^{{\rm u}\ast}$, for simplicity.

We first show the following lemma, which claims that the probability of the alarm flag is invariant as long as the policy conditioned by $\bff_t=0$ is invariant.
\begin{lem}\label{lem:2}
Let $\hpi$ and $\hpi'$ be Markov policies for the augmented system.
If $\hpi_t(\cdot|x_t,0)=\hpi'_t(\cdot|x_t,0)$ for any $t=0,\ldots,T-1$ and $x_t\in\mc{X}$,
then $\Prob^{\hpi}(\bff_t=0)=\Prob^{\hpi'}(\bff_t=0)$ for any $t\in\mc{T}$.
\end{lem}
\begin{IEEEproof}
Since $P^{\hpi}_{\bff}(\bff_t=0|x_{0:t})=0$ for $x_{0:t}\not\in(\Xa^{\rm c})^{t+1}$,
we have
\[
\begin{array}{l}
 \Prob^{\hpi}(\bff_t=0)\\
 = \int_{\mc{X}_{\rm a}^{\rm c}}\int_{\mcA}\cdots\int_{\mc{X}_{\rm a}^{\rm c}} \int_{\mcA} p_t(dx_t|x_{t-1},a_{t-1})\hpi(da_{t-1}|x_{t-1},0)\\
 \quad \cdots p_1(dx_1|x_0,a_0)\hpi(da_0|x_0,0)p_0(dx_0)\\
 = \int_{\mc{X}_{\rm a}^{\rm c}}\int_{\mcA}\cdots\int_{\mc{X}_{\rm a}^{\rm c}} \int_{\mcA} p_t(dx_t|x_{t-1},a_{t-1})\hpi'(da_{t-1}|x_{t-1},0)\\
 \quad \cdots p_1(dx_1|x_0,a_0)\hpi'(da_0|x_0,0)p_0(dx_0)\\
 =\Prob^{\hpi'}(\bff_t=0).
\end{array}
\]
\end{IEEEproof}

Based on Lemma~\ref{lem:2}, we can show the following proposition, which partially characterizes the optimal policy for~\eqref{eq:prob_r2}.
\begin{prop}\label{prop:1}
Let $\hpi^{\ast}$ be the optimal Markov policy for the problem~\eqref{eq:prob_r2}.
Then
\[
 \hpi^{\ast}_t(\cdot|x_t,1)=\pi^{{\rm u}\ast}_t(\cdot|x_t),\quad \forall x_t\in\mcX
\]
\end{prop}
for $t=0,\ldots,T-1$.
\begin{IEEEproof}
For a fixed Markov policy $\hpi\in\hPim$, take $\hpi'\in\hPim$ such that
\[
 \left\{
 \begin{array}{ll}
 \hpi'_t(da_t|x_t,0)\hs :=\hpi_t(da_t|x_t,0),\\
 \hpi'_t(da_t|x_t,1)\hs :=\pi_t^{{\rm u}\ast}(da_t|x_t)
 \end{array}
 \right.
\]
for $t=0,\ldots,T-1$.
Note that $\hpi'$ is feasible for the problem~\eqref{eq:prob_r2} if $\hpi$ is feasible from Lemma~\ref{lem:2}.

Define the value functions associated with $\hpi\in\hPim$ recursively by $V^{\hpi}_T(x_T,\bff_T):=g_T(x_T)$ and
\[
\begin{array}{cl}
 V^{\hpi}_t(x_t,\bff_t) \hs := \int_{\mcA}\{ g_t(x_t,a_t) + \int_{\mcX}\sum_{\bff_{t+1}\in\mc{F}} V^{\hpi}_{t+1}(x_{t+1},\bff_{t+1})\\
  \hs \quad \quad P_{\bff}(\bff_{t+1}|x_t,\bff_t)p_t(dx_{t+1}|x_t,a_t)\} \hpi_t(da_t|x_t,\bff_t)
\end{array}
\]
for $t=0,\ldots,T-1$.
We show that
\begin{equation}\label{eq:valuefunc_comparison}
 V^{\hpi}_t(x_t,\bff_t)\leq V^{\hpi'}_t(x_t,\bff_t),\quad \forall x_t\in\mc{X},\bff_t\in\mc{F}, 
\end{equation}
for any $t\in\mc{T}$ by induction.
It is clear that~\eqref{eq:valuefunc_comparison} holds for $t=T$.
Assume that~\eqref{eq:valuefunc_comparison} holds for $t+1$.
Consider the case with $\bff_t=0$.
Then replacing $\hpi$ with $\hpi'$ yields
\[
 \begin{array}{cl}
 V^{\hpi}_t(x_t,0)=\hs
 \int_{\mcA}\{ g_t(x_t,a_t) + \int_{\mcX}\sum_{\bff_{t+1}\in\mc{F}} V^{\hpi}_{t+1}(x_{t+1},\bff_{t+1})\\
 \hs \quad P_{\bff}(\bff_{t+1}|x_t,\bff_t)p_t(dx_{t+1}|x_t,a_t)\} \hpi'_t(da_t|x_t,0).
\end{array}
\]
From the monotonicity of the Bellman operator, the hypothesis derives
\[
 \begin{array}{cl}
 V^{\hpi}_t(x_t,0)\hs \leq
 \int_{\mcA}\{ g_t(x_t,a_t) + \int_{\mcX}\sum_{\bff_{t+1}\in\mc{F}} V^{\hpi'}_{t+1}(x_{t+1},\bff_{t+1})\\
 \hs \quad P_{\bff}(\bff_{t+1}|x_t,\bff_t)p_t(dx_{t+1}|x_t,a_t)\} \hpi'_t(da_t|x_t,0)\\
 \hs = V^{\hpi'}_t(x_t,0).
\end{array}
\]
On the other hand, for the case with $\bff_t=1$,
\[
 \begin{array}{cl}
 V^{\hpi}_t(x_t,1)=\hs
 \int_{\mcA}\{ g_t(x_t,a_t) + \int_{\mcX} V^{\hpi}_{t+1}(x_{t+1},1)\\
 \hs \quad p_t(dx_{t+1}|x_t,a_t)\} \hpi_t(da_t|x_t,1),
\end{array}
\]
which is the Bellman expectation operator for the unconstrained problem~\eqref{eq:prob_unconstrained}.
Since $\hpi'_t(da_t|x_t,1)$ is the optimal policy for~\eqref{eq:prob_unconstrained}, we get $V^{\hpi}_t(x_t,1)\leq V^{\hpi'}_t(x_t,1)$.

From the inequality~\eqref{eq:valuefunc_comparison}, we have
\[
 J(\hpi)=\Exp[V^{\hpi}_0(x_0,\bff_0)]\leq\Exp[V^{\hpi'}_0(x_0,\bff_0)]=J(\hpi').
\]
Since $\hpi\in\hPih$ can be taken arbitrarily, the claim holds.
\end{IEEEproof}

Proposition~\ref{prop:1} implies that the adversary cares about being detected when there have been no alarms so far, but does no longer care once an alarm has been triggered.
In reality, however, a single alarm may not result in counteractions by the defender due to the presence of false alarms, and a different strategy that avoids serial alarms can possibly be more reasonable.
Therefore, it is more preferable to extend our problem formulation~\eqref{eq:prob_ori} to be able to handle multiple alarms.

\subsection{Multi-alarm Avoidance Policy}

We define the event that alarms are triggered more than or equal to $i$ times,
\[
 \mc{E}^i_{\rm a}:=\{(x_{0:T},a_{0:T-1}): |\mc{T}_{\rm a}(x_{0:T})|\geq i \}
\]
where
\[
 \mc{T}_{\rm a}(x_{0:T}):=\{t\in\mc{T}:x_t\in\Xa\}.
\]

The extended version of the attack impact evaluation problem for multi-alarm avoidance strategies is formulated as follows.
\begin{prob}\label{prob:multi}
The attack impact evaluation problem for multi-alarm avoidance strategies is given by
\begin{equation}\label{eq:prob_multi}
\begin{array}{cl}
 \displaystyle{\max_{\pi\in\Pih}} & J(\pi)\\
 {\rm s.t.} & \Prob^{\pi}(\mc{E}^i_{\rm a})\leq \Delta_i,\quad i=1,\ldots,M
\end{array}
\end{equation}
with given constants $\Delta_i\geq 0$ for $i=1,\ldots,M$.
\end{prob}

The same idea of the alarm flag state augmentation proposed in Sec.~\ref{sec:approach} can be applied to Problem~\ref{prob:multi} as well by augmenting information on the number of alarms instead of the binary information.
The augmented state space and the augmented system for Problem~\ref{prob:multi} are defined as follows.
\begin{defin}
The augmented state space of $\mcX$ for Problem~\ref{prob:multi} is defined as $\hat{\mcX}:=\mcX\times \mc{F}$
with $\mc{F}:=\{0,\ldots,M\}$.
The augmented system is defined as
\[
\left\{
\begin{array}{cl}
 x_{t+1} \hs = f(x_t,a_t,w_t),\\
 \bff_{t+1} \hs = \left\{
 \begin{array}{ll}
 M & {\rm if}\ \bff_t=M,\\
 \bff_t+1 & {\rm if}\ x_{t+1}\in\Xa\ {\rm and}\ \bff_t<M,\\
 \bff_t & {\rm otherwise}.
 \end{array}
 \right.
\end{array}
\right.
\]
\end{defin} 

The alarm number augmentation naturally leads to an equivalent problem
\begin{equation}\label{eq:prob_rr1}
\begin{array}{cl}
\displaystyle{\max_{\hpi\in\hPim}} & J(\hpi)\\
 {\rm s.t.} & \Prob^{\hpi}(\bff_T\geq i) \leq \Delta_i,\quad i=1,\ldots,M
\end{array}
\end{equation}
where the search space is the set of Markov policies for the state space augmented with the number of alarms.
The following theorem is the correspondence of Theorem~\ref{thm:1}.

\begin{theorem}\label{thm:2}
The optimal values of the problems~\eqref{eq:prob_multi} and~\eqref{eq:prob_rr1} are equal.
\end{theorem}
\begin{IEEEproof}
The claim can be proven in a manner similar to that of Theorem~\ref{thm:1}.
\end{IEEEproof}

Moreover, the correspondence of Proposition~\ref{prop:1} is described as follows.
\begin{prop}\label{prop:2}
Let $\hpi^{\ast}$ be the optimal Markov policy for the reformulated problem~\eqref{eq:prob_rr1}.
Then
\[
 \hpi^{\ast}_t(da_t|x_t,M)=\pi^{{\rm u}\ast}_t(da_t|x_t),\quad \forall x_t\in\mcX
\]
for $t=0,\ldots,T-1$.
\end{prop}
\begin{IEEEproof}
The claim can be proven in a manner similar to that of Proposition~\ref{prop:1}.
\end{IEEEproof}
Proposition~\ref{prop:2} means that the adversary does not avoid detection after the number of alarms reaches $M$.

{\it Remark:}
The constraint in the extended problem~\eqref{eq:prob_multi} restricts the probability distribution of the number of alarms.
In other words, the formulation utilizes a risk measure on a probability distribution instead of a typical statistic.
Several risk measures have been proposed, such as CVaR, which is one of the most commonly used coherent risk measures~\cite{Artzner1999Coherent}.
Those risk measures compress risk of a random variable with a distribution into a scalar value.
Because our formulation uses the full information of the distribution, the constraint can be regarded as a fine-grained version of standard risk measures.

\section{NUMERICAL EXAMPLE}
\label{sec:num}


\subsection{Simulation Setup}
Consider the one-dimensional discrete-time integrator
\[
 z_{t+1}=z_t+a_t+w_t,\ z_0=0,
\]
with the CUSUM attack detector~\cite[Chap.~2]{Basseville1993Detection}
\[
\left\{
\begin{array}{cl}
 z^{\rm d}_{t+1} \hs = \max(0,z^{\rm d}_t+|x_t|-b^{\rm d}),\\
 \delta_t \hs = \left\{
 \begin{array}{cl}
 1 \hs {\rm if}\ z^{\rm d}_t\geq\tau^{\rm d},\\
 0 \hs {\rm otherwise},
 \end{array}
 \right.
\end{array}\quad z^{\rm d}_0=0,
\right.
\]
with the bias $b^{\rm d}>0$ and the threshold $\tau^{\rm d}>0$, where $\mcA=\mc{W}=\mathbb{R}$.
The state space $\mc{X}=\mathbb{R}^2$ and the alarm region $\Xa$ are constructed according to Sec.~\ref{subsec:sys}.
The process noise follows the white Gaussian distribution with mean zero and variance $\sigma^2$.
The adversary's objective is to drive the system state around a reference value $z_{\rm ref}\in\mathbb{R}$.
Accordingly, the objective function is set to a quadratic function
\[
 \textstyle{J(\pi):=-\Exp^\pi\left[ \sum_{t=0}^T (z_t-z_{\rm ref})^2 \right].}
\]
The constants are specifically set to
$T=15,$ $\sigma=0.1,$ $b^{\rm d}=0.8,$ $\tau^{\rm d}=2,$ and $z_{\rm ref}=1.5$.
We compute $\hpi^{\ast}_t(da_t|x_t,0)$ based on discretization of the state and input spaces and use a standard linear programming approach for solving the resulting constrained finite MDP~\cite{Altman1999Constrained}. 
On the other hand, we analytically compute $\hpi^{\ast}_t(da_t|x_t,1)$ as the unconstrained discrete-time linear quadratic regulator~\cite[Chap.~4]{Bertsekas2012Dynamic} based on Proposition~\ref{prop:1}.

\subsection{Simulation Results}
\begin{figure}[t]
  \centering
  \subfloat[][]{
    \includegraphics[width=.98\linewidth]{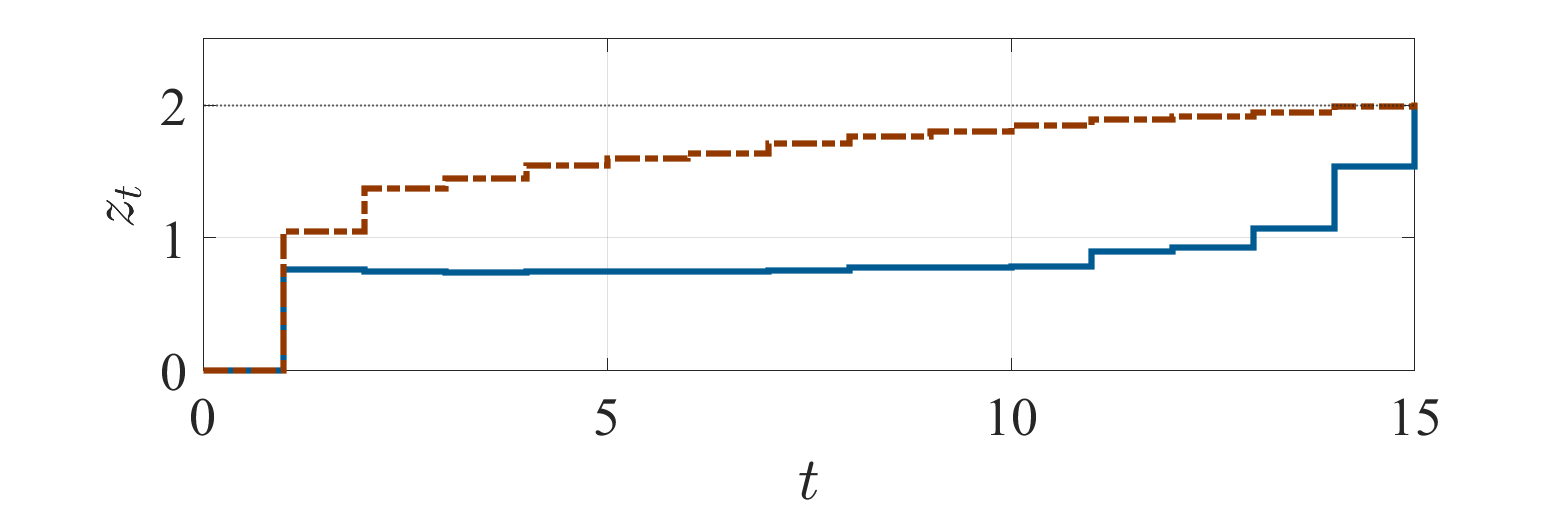}\label{subfig:single_zt}
    }\\
  \subfloat[][]{
    \includegraphics[width=.98\linewidth]{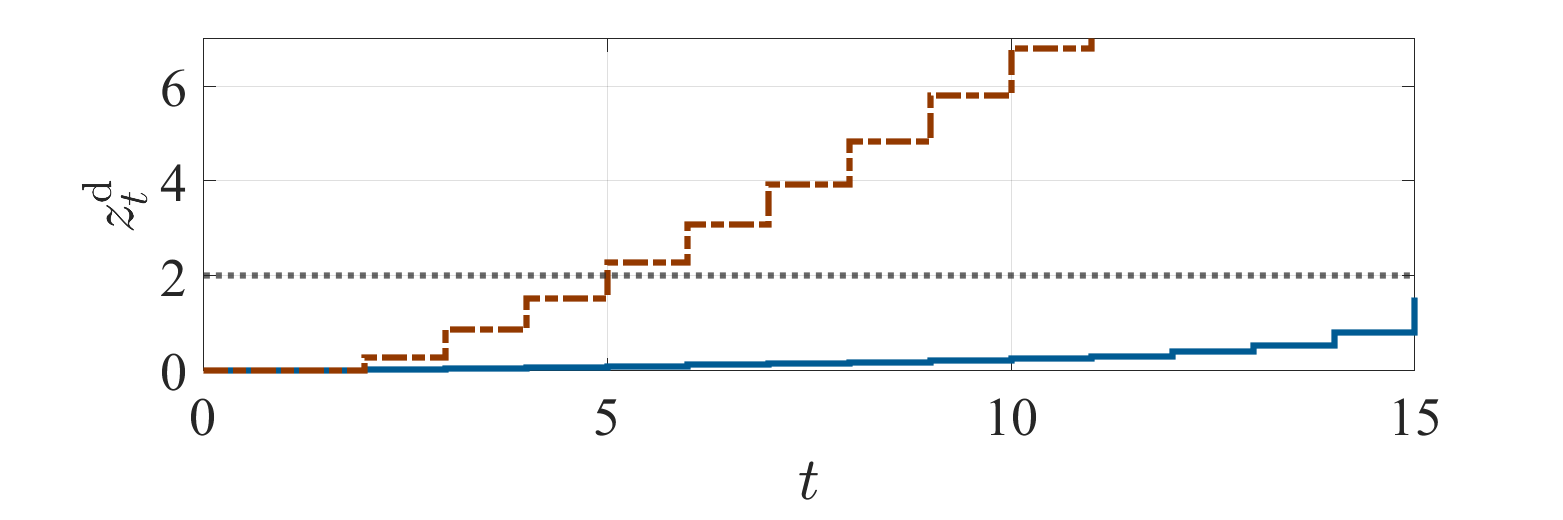}\label{subfig:single_zdt}
    }\\
  \subfloat[][]{
    \includegraphics[width=.98\linewidth]{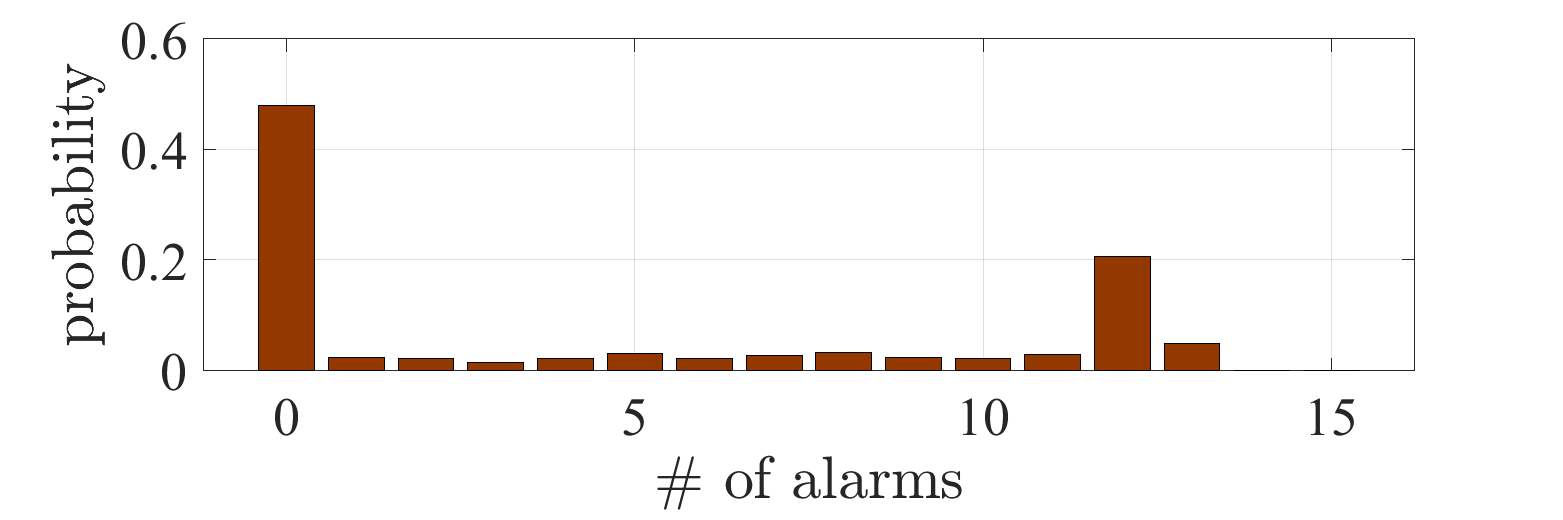}\label{subfig:single_number}
    }
    \caption{
    Simulation results with the formulation of Problem~\ref{prob:single}.
    (a): Conditional means of $z_t$. The dash line corresponds to the case where the alarm is triggered at least once, while the solid line corresponds to the case where the alarm is not triggered during the process.
    The dotted line describes the reference value $z_{\rm ref}$.
    (b): Conditional means of $z^{\rm d}_t$. The dotted line describes the threshold $\tau^{\rm d}$.
    (c): Probabilities with respect to the total number of alarms during the process.
    }
    \label{fig:single}
\end{figure}

We first treat the formulation of Problem~\ref{prob:single}.
Set the constant on the stealth condition by $\Delta=0.5$.
Fig.~\ref{fig:single} shows the simulation results with the optimal policy obtained by solving the equivalent problem~\eqref{eq:prob_r2}.
Figs.~\ref{subfig:single_zt} and~\ref{subfig:single_zdt} depict the empirical means of $z_t$ and $z^{\rm d}_t$ conditioned by whether an alarm is triggered at least once during the process, or not, respectively.
It can be observed that $z^{\rm d}_t$ increases even after an alarm is triggered, as claimed in Sec.~\ref{subsec:opt_pol_after_alarm}.
Fig.~\ref{subfig:single_number} depicts the probabilities with respect to the total number of alarms during the process.
It can be observed that a large number of alarms occur with a high probability.
The result indicates that the formulation of Problem~\ref{prob:single} leads to a policy such that the number of alarms becomes large once an alarm is triggered.

\begin{figure}[t]
  \centering
  \subfloat[][]{
    \includegraphics[width=.98\linewidth]{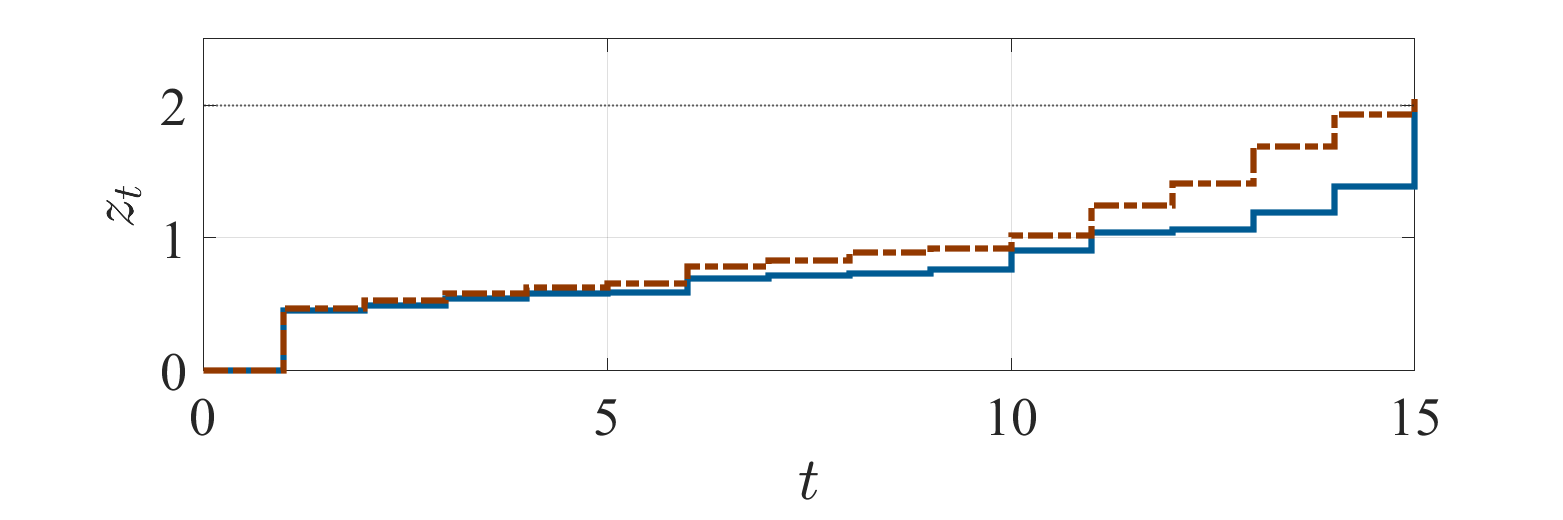}\label{subfig:multi_zt}
    }\\
  \subfloat[][]{
    \includegraphics[width=.98\linewidth]{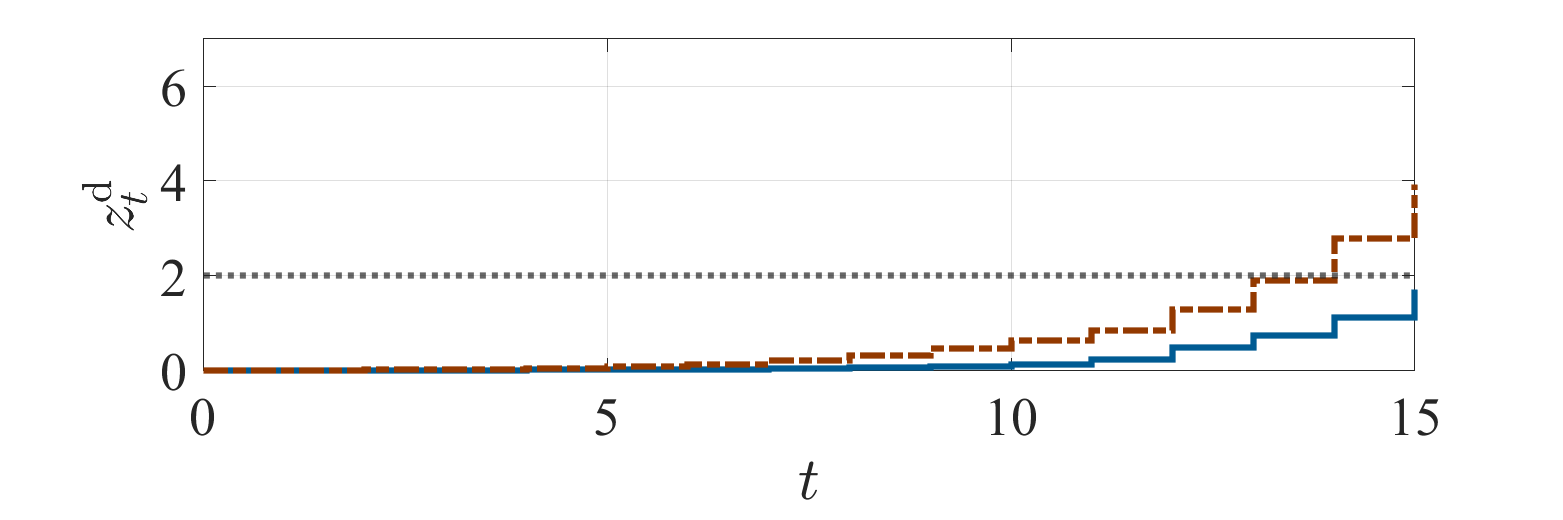}\label{subfig:multi_zdt}
    }\\
  \subfloat[][]{
    \includegraphics[width=.98\linewidth]{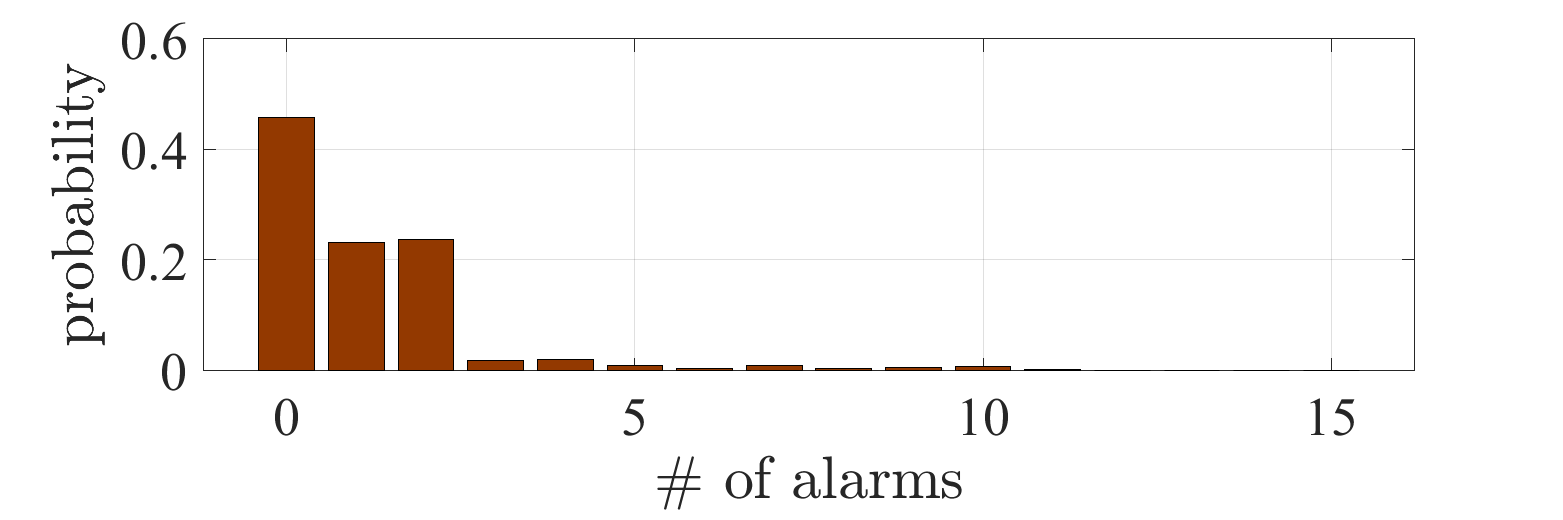}\label{subfig:multi_number}
    }
    \caption{
    Simulation results with the formulation of Problem~\ref{prob:multi}.
    Each subfigure corresponds to those in Fig.~\ref{fig:single}.
    }
    \label{fig:multi}
\end{figure}

We next treat the formulation of Problem~\ref{prob:multi}.
Set the constants on the stealth condition by
\[
 \Delta_1=0.5,\quad \Delta_2=0.3,\quad \Delta_3=0.1
\]
with $M=3$.
Fig.~\ref{fig:multi} shows the simulation results.
The subfigures correspond to those in Fig.~\ref{fig:single}.
It can be observed that the trajectory of $z^{\rm d}_t$ is kept less than the detection threshold $\tau^{\rm d}$ over almost the entire period.
Accordingly, the probability depicted in Fig.~\ref{fig:multi} suggests that the obtained policy avoids multiple alarms.

\section{CONCLUSION}
\label{sec:conc}

This note has addressed the attack impact evaluation problem for stochastic control systems.
The problem is formulated as an optimal control problem with a temporally joint chance constraint.
The difficulty to solve the optimal control problem lies in the explosion of the search space owing to the dependency of the optimal policy on the entire history.
In this note, we have shown that the information whether an alarm has been triggered or not is sufficient for determining the optimal decision.
By augmenting the alarm flag with the original state space, we can obtain an equivalent optimal control problem in a computationally tractable form.
Moreover, the formulation is extended to handle multi-alarm avoidance policies by taking the number of alarms into account.

Future research directions include development of a numerical algorithm that efficiently solves the reformulated problem.
Although the search space is hugely reduced by the proposed method, it still suffers from the curse of dimensionality occurring from space discretization.
In addition, it is interesting to clarify the relationship between the chance constraint considered in our formulation and existing risk measures, such as CVaR.
Although we have used full information of the probability distribution, coherent risk measures can effectively compress the information, the property of which can possibly be used for efficient numerical algorithms.

\begin{figure}[t]
  \centering
  \includegraphics[width=0.98\linewidth]{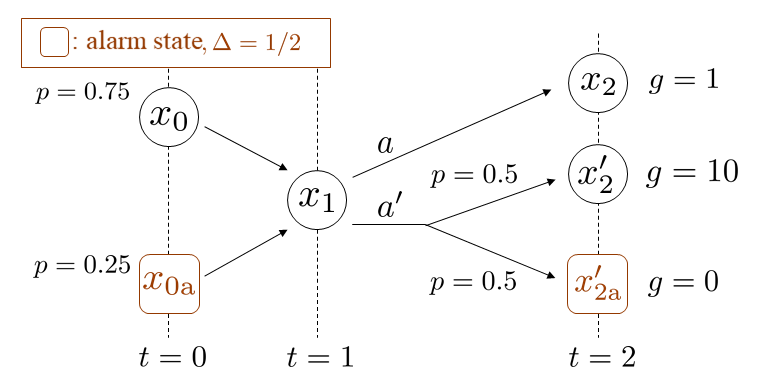}
  \caption{
  MDP for which Markov policies cannot attain the optimal value of Problem~\ref{prob:single}.
  The adversary can inject an input only at $t=1$.
  The attack impact is determined by the state at the final time step.
  }
  \label{fig:ex_MDP}
\end{figure}

\appendix
We provide an example for which Markov policies cannot attain the optimal value in the formulation~\eqref{eq:prob_ori}.
Consider the finite MDP illustrated by Fig.~\ref{fig:ex_MDP}.
The adversary can inject an input only at $t=1$.
When the input $a$ is selected, the state reaches $x_2$ with probability one and the resulting attack impact is $1$.
On the other hand, when the input $a'$ is selected, the state reaches $x'_2$ or $x'_{2{\rm a}}$ with equal probabilities.
The impact is 10 in the case of $x'_2$, while there is no impact in the case of $x'_{2{\rm a}}$.
The alarm region is given as $\mc{X}_{\rm a}=\{x_{0{\rm a}},x'_{2{\rm a}}\}$.
The input $a'$ can be interpreted as a risky action in the sense that it leads to large impact in expectation but may trigger an alarm.

We derive the optimal policy.
The history-dependent policies can be parameterized by $\pi(a'|x_0,x_1)=\alpha$ and $\pi(a'|x_{0{\rm a}},x_1)=\beta$ with parameters $(\alpha,\beta)\in[0,1]\times[0,1]$.
The joint chance constraint is written by
$\Prob^{\pi}(\lor_{t\in\mc{T}}X_t\in\mc{X}_{\rm a})\leq \Delta \Leftrightarrow \Prob^{\pi}(x_{0{\rm a}}) + \Prob^{\pi}(x_1)\pi(a'|x_0,x_1)p(x'_{2{\rm a}}|x_1,a')\leq 1/2\Leftrightarrow 1/4+3/4 \cdot \alpha/2\leq 1/2\Leftrightarrow \alpha\leq 2/3.$
Thus the feasible region of $(\alpha,\beta)$ is $[0,2/3]\times[0,1]$.
The objective function is written by
$J(\pi)=\Prob^{\pi}(x_2)+\Prob^{\pi}(x'_2)\cdot 10= 3/4(1-\alpha) + (1-\beta)/4+(3/8\cdot \alpha+\beta/8)\cdot10= 3\alpha+\beta+1.$
Because this is monotonically increasing with respect to $\alpha$ and $\beta$, the optimal values are $(\alpha^{\ast},\beta^{\ast})=(2/3,1),$ which leads to
\[
 \left\{
\begin{array}{ll}
 \pi^{\ast}(a'|x_0,x_1)\hs =2/3,\\
 \pi^{\ast}(a'|x_{0{\rm a}},x_1)\hs =1.
\end{array}
\right.
\]
On the other hand, the Markov policies can be parameterized by $\pi^{\rm m}(a'|x_1)=\gamma$ with $\gamma\in[0,1]$.
The joint chance constraint imposes $\gamma\leq 2/3$ and the objective function is $4\gamma+1$.
Thus the optimal policy is given by $\gamma^{\ast}=2/3,$ namely $\pi^{{\rm m}\ast}(a'|x_1)=2/3.$
Denoting the value of the objective function with a policy $\pi$ by $J(\pi),$ we have $J(\pi^{\ast})=4>11/3=J(\pi^{{\rm m}\ast}),$ which implies that Markov policies cannot attain the optimal value for this instance.

The optimal history-dependent policy means that the adversary reduces the risk when no alarm has been triggered while she selects the risky input when an alarm has been triggered.
In other words, the decision making at the state $x_1$ depends on the alarm flag.
This observation leads to the hypothesis that this binary information in addition to the current state is sufficient for optimal decision making, giving rise to the idea of alarm flag state augmentation.

\bibliography{sshrrefs}

\begin{thebibliography}{10}
\providecommand{\url}[1]{#1}
\csname url@samestyle\endcsname
\providecommand{\newblock}{\relax}
\providecommand{\bibinfo}[2]{#2}
\providecommand{\BIBentrySTDinterwordspacing}{\spaceskip=0pt\relax}
\providecommand{\BIBentryALTinterwordstretchfactor}{4}
\providecommand{\BIBentryALTinterwordspacing}{\spaceskip=\fontdimen2\font plus
\BIBentryALTinterwordstretchfactor\fontdimen3\font minus
  \fontdimen4\font\relax}
\providecommand{\BIBforeignlanguage}[2]{{%
\expandafter\ifx\csname l@#1\endcsname\relax
\typeout{** WARNING: IEEEtran.bst: No hyphenation pattern has been}%
\typeout{** loaded for the language `#1'. Using the pattern for}%
\typeout{** the default language instead.}%
\else
\language=\csname l@#1\endcsname
\fi
#2}}
\providecommand{\BIBdecl}{\relax}
\BIBdecl

\bibitem{History2018Hemsley}
K.~E. Hemsley and D.~R. E.~Fisher, ``History of industrial control system cyber
  incidents,'' U.S. Department of Energy Office of Scientific and Technical
  Information, Tech. Rep. INL/CON-18-44411-Rev002, 2018.

\bibitem{Nicolas2011Stuxnet}
N.~Falliere, L.~O. Murchu, and E.~Chien, ``W32. {Stuxnet Dossier},'' Symantec,
  Tech. Rep., 2011.

\bibitem{CISA2014}
{Cybersecurity \& Infrastructure Security Agency}, ``Stuxnet malware
  mitigation,'' Tech. Rep. ICSA-10-238-01B, 2014, [Online]. Available:
  \url{https://www.us-cert.gov/ics/advisories/ICSA-10-238-01B}.

\bibitem{CISA2017}
------, ``{HatMan} - safety system targeted malware,'' Tech. Rep.
  MAR-17-352-01, 2017, [Online]. Available:
  \url{https://www.us-cert.gov/ics/MAR-17-352-01-HatMan-Safety-System-Targeted-Malware-Update-B}.

\bibitem{CISA2018}
------, ``Cyber-attack against {U}krainian critical infrastructure,'' Tech.
  Rep. IR-ALERT-H-16-056-01, 2018, [Online]. Available:
  \url{https://www.us-cert.gov/ics/alerts/IR-ALERT-H-16-056-01}.

\bibitem{Kaplan1981On}
S.~Kaplan and B.~J. Garrick, ``On the quantitative definition of risk,''
  \emph{Risk Analysis}, vol.~1, no.~1, pp. 11--27, 1981.

\bibitem{Sridhar2012Cyber}
S.~{Sridhar}, A.~{Hahn}, and M.~{Govindarasu}, ``Cyber–physical system
  security for the electric power grid,'' \emph{Proc. IEEE}, vol. 100, no.~1,
  pp. 210--224, 2012.

\bibitem{Teixeira2015Secure}
A.~{Teixeira}, K.~C. {Sou}, H.~{Sandberg}, and K.~H. {Johansson}, ``Secure
  control systems: A quantitative risk management approach,'' \emph{IEEE
  Control Systems Magazine}, vol.~35, no.~1, pp. 24--45, Feb 2015.

\bibitem{Murguia2018Reachable}
C.~Murguia and J.~Ruths, ``On reachable sets of hidden {CPS} sensor attacks,''
  in \emph{2018 Annual American Control Conference (ACC)}, 2018, pp. 178--184.

\bibitem{Murguia2019On}
C.~{Murguia} and J.~{Ruths}, ``On model-based detectors for linear
  time-invariant stochastic systems under sensor attacks,'' \emph{IET Control
  Theory Applications}, vol.~13, no.~8, pp. 1051--1061, 2019.

\bibitem{Alvaro2011Attacks}
A.~A. C\'{a}rdenas, S.~Amin, Z.-S. Lin, Y.-L. Huang, C.-Y. Huang, and
  S.~Sastry, ``Attacks against process control systems: {R}isk assessment,
  detection, and response,'' in \emph{Proc. the 6th ACM ASIA Conference on
  Computer and Communications Security}, 2011.

\bibitem{Bai2017Data}
C.-Z. Bai, F.~Pasqualetti, and V.~Gupta, ``Data-injection attacks in stochastic
  control systems: {D}etectability and performance tradeoffs,''
  \emph{Automatica}, vol.~82, pp. 251 -- 260, 2017.

\bibitem{Mo2016On}
Y.~{Mo} and B.~{Sinopoli}, ``On the performance degradation of cyber-physical
  systems under stealthy integrity attacks,'' \emph{IEEE Trans. Autom.
  Control}, vol.~61, no.~9, pp. 2618--2624, Sep. 2016.

\bibitem{Umsonst2017Security}
D.~{Umsonst}, H.~{Sandberg}, and A.~A. {C\'{a}rdenas}, ``Security analysis of
  control system anomaly detectors,'' in \emph{Proc. 2017 American Control
  Conference (ACC)}, May 2017, pp. 5500--5506.

\bibitem{Hirzallah2018Computation}
N.~H. Hirzallah and P.~G. Voulgaris, ``On the computation of worst attacks: {A}
  {LP} framework,'' in \emph{2018 Annual American Control Conference (ACC)},
  2018, pp. 4527--4532.

\bibitem{Chen2018Optimal}
Y.~{Chen}, S.~{Kar}, and J.~M.~F. {Moura}, ``Optimal attack strategies subject
  to detection constraints against cyber-physical systems,'' \emph{IEEE Trans.
  Contr. Netw. Systems}, vol.~5, no.~3, pp. 1157--1168, 2018.

\bibitem{Teixeira2019Optimal}
A.~M.~H. Teixeira, ``Optimal stealthy attacks on actuators for strictly proper
  systems,'' in \emph{2019 IEEE 58th Conference on Decision and Control (CDC)},
  2019, pp. 4385--4390.

\bibitem{Milosevic2019Estimating}
J.~Milo{\v{s}}evi{\'c}, H.~Sandberg, and K.~H. Johansson, ``Estimating the
  impact of cyber-attack strategies for stochastic networked control systems,''
  \emph{IEEE Trans. Control Netw. Syst.}, vol.~7, no.~2, pp. 747--757, 2019.

\bibitem{Fang2019Stealthy}
C.~Fang, Y.~Qi, J.~Chen, R.~Tan, and W.~X. Zheng, ``Stealthy actuator signal
  attacks in stochastic control systems: {P}erformance and limitations,''
  \emph{IEEE Trans. Autom. Control}, vol.~65, no.~9, pp. 3927--3934, 2019.

\bibitem{Sui2020The}
T.~Sui, Y.~Mo, D.~Marelli, X.~Sun, and M.~Fu, ``The vulnerability of
  cyber-physical system under stealthy attacks,'' \emph{IEEE Trans. Autom.
  Control}, vol.~66, no.~2, pp. 637--650, 2020.

\bibitem{Wang2022Optimal}
X.-L. Wang, G.-H. Yang, and D.~Zhang, ``Optimal stealth attack strategy design
  for linear cyber-physical systems,'' \emph{IEEE Trans. on Cybern.}, vol.~52,
  no.~1, 2022.

\bibitem{Khazraei2022Resiliency}
A.~Khazraei, H.~Pfister, and M.~Pajic, ``Resiliency of perception-based
  controllers against attacks,'' in \emph{Proc. Learning for Dynamics and
  Control Conference}, 2022, pp. 713--725.

\bibitem{Bauerle2011Markov}
N.~B{\"a}uerle and J.~Ott, ``Markov decision processes with
  average-value-at-risk criteria,'' \emph{Mathematical Methods of Operations
  Research}, vol.~74, no.~3, pp. 361--379, 2011.

\bibitem{Haskell2015Convex}
W.~B. Haskell and R.~Jain, ``A convex analytic approach to risk-aware {M}arkov
  decision processes,'' \emph{SIAM Journal on Control and Optimization},
  vol.~53, no.~3, pp. 1569--1598, 2015.

\bibitem{Chow2017Risk}
Y.~Chow, M.~Ghavamzadeh, L.~Janson, and M.~Pavone, ``Risk-constrained
  reinforcement learning with percentile risk criteria,'' \emph{The Journal of
  Machine Learning Research}, vol.~18, no.~1, pp. 6070--6120, 2017.

\bibitem{Baier2008Principles}
C.~Baier and J.-P. Katoen, \emph{Principles of Model Checking}.\hskip 1em plus
  0.5em minus 0.4em\relax MIT Press, 2008.

\bibitem{Ono2012Joint}
M.~Ono, Y.~Kuwata, and J.~Balaram, ``Joint chance-constrained dynamic
  programming,'' in \emph{Proc. IEEE Conference on Decision and Control (CDC)},
  2012, pp. 1915--1922.

\bibitem{Ono2015Chance}
M.~Ono, M.~Pavone, Y.~Kuwata, and J.~Balaram, ``Chance-constrained dynamic
  programming with application to risk-aware robotic space exploration,''
  \emph{Autonomous Robots}, vol.~39, no.~4, pp. 555--571, 2015.

\bibitem{Thorpe2022Data}
A.~Thorpe, T.~Lew, M.~Oishi, and M.~Pavone, ``Data-driven chance constrained
  control using kernel distribution embeddings,'' in \emph{Proc. Learning for
  Dynamics and Control Conference}, 2022, pp. 790--802.

\bibitem{Patil2022ChanceA}
A.~Patil, A.~Duarte, A.~Smith, F.~Bisetti, and T.~Tanaka, ``Chance-constrained
  stochastic optimal control via path integral and finite difference methods,''
  in \emph{Proc. IEEE Conference on Decision and Control (CDC)}, 2022, pp.
  3598--3604.

\bibitem{Patil2022ChanceB}
A.~Patil, A.~Duarte, F.~Bisetti, and T.~Tanaka, ``Chance-constrained stochastic
  optimal control via {HJB} equation with {D}irichlet boundary condition,''
  2022, [Online]. Available:
  \url{https://sites.utexas.edu/tanaka/files/2022/07/Chance_Constrained_SOC.pdf}.

\bibitem{Sasahara2022Attack}
H.~Sasahara, T.~Tanaka, and H.~Sandberg, ``Attack impact evaluation by exact
  convexification through state space augmentation,'' in \emph{Proc. IEEE
  Conference on Decision and Control (CDC)}, 2022, pp. 7084--7089.

\bibitem{Hinderer1970Foundations}
K.~Hinderer, \emph{Foundations of Non-stationary Dynamic Programming with
  Discrete Time Parameter}.\hskip 1em plus 0.5em minus 0.4em\relax Springer,
  1970.

\bibitem{Bertsekas1996Stochastic}
D.~Bertsekas and S.~Shreve, \emph{Stochastic Optimal Control: The Discrete-Time
  Case}.\hskip 1em plus 0.5em minus 0.4em\relax Athena Scientific, 1996.

\bibitem{Munos2002Variable}
R.~Munos and A.~Moore, ``Variable resolution discretization in optimal
  control,'' \emph{Machine Learning}, vol.~49, no.~2, pp. 291--323, 2002.

\bibitem{Artzner1999Coherent}
P.~Artzner, F.~Delbaen, J.-M. Eber, and D.~Heath, ``Coherent measures of
  risk,'' \emph{Mathematical Finance}, vol.~9, no.~3, pp. 203--228, 1999.

\bibitem{Basseville1993Detection}
M.~Basseville, I.~V. Nikiforov \emph{et~al.}, \emph{Detection of Abrupt
  Changes: Theory and Application}.\hskip 1em plus 0.5em minus 0.4em\relax
  Prentice Hall Englewood Cliffs, 1993.

\bibitem{Altman1999Constrained}
E.~Altman, \emph{Constrained Markov Decision Processes}.\hskip 1em plus 0.5em
  minus 0.4em\relax Chapman and Hall/CRC, 1999.

\bibitem{Bertsekas2012Dynamic}
D.~Bertsekas, \emph{Dynamic Programming and Optimal Control: Volume I}.\hskip
  1em plus 0.5em minus 0.4em\relax Athena Scientific, 2012.

\end{thebibliography}
\bibliographystyle{IEEEtran}

\end{document}